\documentclass{article}
\usepackage[utf8]{inputenc}
\usepackage{comment}
\usepackage{changepage}
\usepackage{multirow}
\usepackage[citecolor=blue,urlcolor=blue,bookmarks=false,hypertexnames=true]{hyperref} 
\usepackage{graphicx}
\usepackage{amsmath}
\title{Optimal Strategies to Catch Randomly Walking Cat}
\author{Rüdiger Jehn}

\begin{document}

\maketitle

{\bf Abstract}\\

The optimal strategies to catch a randomly walking cat in various environments are presented. All games have a player that opens a box at step $i$. If the cat is in this box the player wins, if not, the cat moves randomly to an adjacent box, or in case the box is open to the outside, the cat may escape. If the cat is not escaping, the next step of the game is played. In case the cat has doors to escape, the optimal strategy is determined that minimizes the escape chances of the cat. If there are no doors to the outside, the strategies are calculated that minimize the game duration.\\

The environments are 1) a one-dimensional array of up to 9 boxes in a line 2) this line is connected to a ring and 3) a $2 \times m$ grid (with $2\le m\le 4$) of boxes. In cases 1) and 3) the boxes may or may not have exits to the outside. Numerical proofs for the optimality of the presented strategies are outlined.\\

Also a formula is presented for the average number of steps it takes a randomly walking cat to exit a $2\times m$ grid when no player is involved. This formula is only based on Fibonacci numbers if $m$ is even and Lucas numbers if $m$ is odd. 

\section{Introduction}

The textbook {\it Principles of Random Walk} by Spitzer~\cite{Spitzer1964}, first published in 1964, provides all the mathematics related to random walk processes. But it does not include cases, where a player is involved that can remove the walking object.\\

\begin{figure}[h]
\centering
\includegraphics[width=1.0\linewidth]{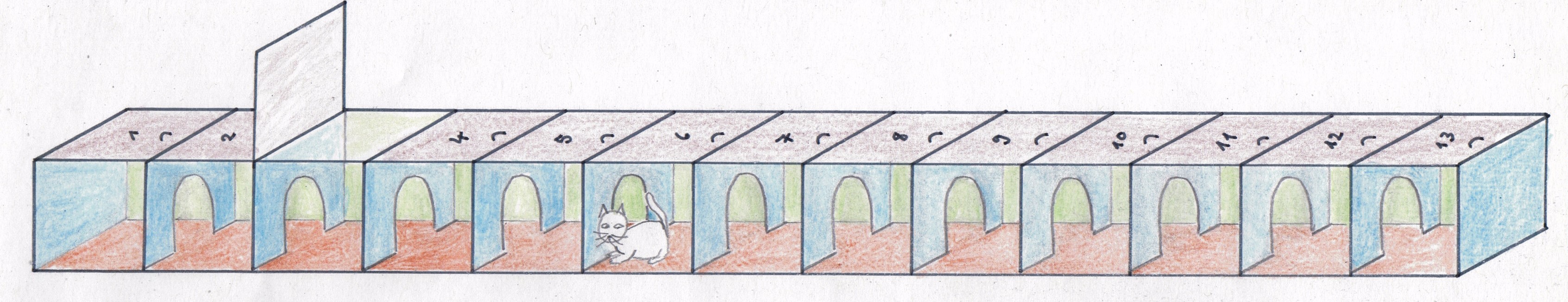}
\label{fig:cat}
\vspace{-6mm}
\caption{At each step a player opens one of $n$ boxes to search for a cat which changes its position in a simple random walk after every non successful opening of a box.}
\end{figure}

In all problems studied here, a cat is hiding in one of $n$ boxes and a single player tries to catch the cat. In one step, first the player opens a box, if the cat is found, the game is won, if not the cat will move randomly to an adjacent box. In case where the boxes have exits to the "outside", the cat may escape. The cat is moving in any of the possible directions with the same probability, i.e.~if there are two options each option will be chosen with a 50~\% probability, if the cat is in box~1 which is only connected with box~2, it will move with 100~\% probability to box~2. The cat is initially assigned a random box. Fig.\ref{fig:cat} shows an example for a typical set-up.

Strategies are denoted as $s_1s_2\ldots s_k\overline{t_1t_2\ldots t_\ell}$, which means the boxes $s_1$ to $s_k$ are opened in the steps from 1 to $k$ and then the boxes $t_1$ to $t_\ell$ are opened repeatedly in this sequence until the cat is found or, in case where the outer boxes have escape doors, the cat has escaped.\\

Since the problems are symmetric, all strategies have one or more mirror strategies. In the case of one-dimensional arrangement of the boxes like in Fig.\ref{fig:cat}, strategy $s_1s_2\ldots s_k$ will achieve the same results as strategy $n-s_1,n-s_2,\ldots n-s_k$. In the case of a two-dimensional grid of boxes the strategies that start with opening a box in a corner have 3 mirror strategies, and a strategy not starting with opening a box in a corner has one mirror strategy. In this paper, the strategies where the first box opened is on the left side ($s_1\le\frac{n}{2}$) are presented.

\section{Walking grid with no exits}

\subsection{Strategy to find cat in a minimum number of steps}

The problem where a cat or a princess changes every night her box or room can be found in many places in the internet. Sequentially opening the boxes 2, 3, 4, ..., $n-1$, $n-1$, ... 4, 3, 2 guarantees that the cat (or princess) is found latest after step $2n-4$ ($n>2$). A nice illustration is given by Perfect~\cite{perfect}.\\

On average with this strategy it takes approximately $n-1.5$ steps to find the cat:
for n odd, the probability to find the cat in the first step is $\frac{1}{n}$, in the next steps up to step $n-2$ the probability is $\frac{1}{2n}$, at step $n-1$ the probability is about $\frac{n+1}{n(n-1)}$ and in the final steps the probability is about $\frac{n+1}{2n(n-1)}$. Therefore, the expectation value $E_{save}$ is
$$E_{save}(n) = \frac{1}{n} \left ( 1 + \frac{1}{2} \sum_{i=2}^{n-2} i + n+1 + \frac{1}{2} \frac{n+1}{n -1}\sum_{i=2}^{n-2} (i + n - 2) \right ).$$
This yields
$$E_{save}(n) = \frac{1}{n}\frac{2n^3-5n^2+3n+2}{2(n-1)}\approx n - 1.5$$

For n even, it takes a bit longer, but eventually $E_{save}(n)$ is also approaching $n-1.5$ ($E_{save}(1000)= 998.47$ vs $E_{save}(1001)= 999.495$).\\

If you open the boxes randomly, even if the search could go on forever, the average number of steps to find the cat is only 1.5 larger:
$$E_{random}(n) = \frac{1}{n} \sum_{i=1}^\infty i \cdot (1-\frac{1}{n})^{i-1}  = n $$

\subsection{Strategy to find cat in fastest way}

The strategy presented in Section~2.1 guarantees to find the cat in no more than $2n-4$ steps, but it is not the fastest strategy if $n>4$. A semi-automated algorithm was written, that tries all possible strategies up to a length of $maxdepth$ steps, but prunes the branches where a given threshold for the average duration of the search is surpassed. Whenever a new minimum is found, the threshold is reduced.\\

\vspace{-1mm}
The proof that a strategy is optimal is also performed with a search tree algorithm. For each step, it is tested whether opening a different box allows to reduce the expected duration of the game. If all alternatives lead to an increase of the game duration, then it is proved that for this step the optimum was found. Since the repetitive pattern of the strategies eventually also lead to a repetitive pattern of the cat residence probability distribution, it is only necessary to prove the optimality up to the step where the repetition starts.\\

\vspace{-1mm}
The optimal strategies that minimize the duration of the search are presented in Table~\ref{catch_fast} for the number of boxes $n$ up to 7. For $n=8$ only the first 87 steps of the optimal strategy are known so far.\\

\begin{table}[h!]
\begin{center}
\begin{tabular}{|c|c|c|c|}
\hline
\textbf{$n$} & \textbf{Strategy} & \textbf{Average duration} & \textbf{Duration of strategy 2.1} \\
\hline
2 & 11 & 3/2 & \\
3 & 22 & 5/3 & 5/3\\
4 & 2332 & 39/16 (2.4375) & 39/16 (2.4375) \\
5 & $\overline{2442}$ & 44/15 (2.9$\overline{3}$) & 71/20 (3.55) \\
6 & 255233$\overline{5522}$ & 34165/9984 (3.4220) & 279/64 (4.3594)\\ 
7 & 263265432$\overline{6325}$ & 7373/1792 (4.11440) & 9897/1792 (5.5229) \\
8 & $s_8$ & 4.74959 & 25963/4096 (6.339) \\
\hline
\end{tabular}
\end{center}
\vspace{-4mm}
\caption{Optimal strategy (or one of the possible optimal strategies) to catch as fast as possible the cat performing a random walk through $n$ boxes, the average game duration of this strategy and the average game duration of the strategy 234...$n$-1$n$-1...432 where the cat is guaranteed to be found within $2n-4$ steps. The 87 known digits of $s_8$ are given below in the text.}
\label{catch_fast}
\end{table}

\vspace{-1mm}
The cases $n$ = 2 and 3 are trivial. For $n=4$ also strategy 22332 gives the same expected duration of 39/16.\\

\vspace{-1mm}
But for $n>4$ alternative strategies are found where the cat is found faster. For $n=5$, strategy $\overline{2442}$ produces a relative cat residence probability of 
\newline
$[0.3, 0.1, 0.3, 0.3, 0]$ in steps 2, 6, 10, ..., hence the repetition cycle has a length of 4 and the proof of optimality only needs to be made for steps 1 to 5. The game situations in step 2 and 6 are identical.\\

\vspace{-1mm}
For $n=6$, strategy 255233$\overline{5522}$ produces a relative cat residence probability of $\frac{1}{66}[0, 20, 13, 15, 13, 5]$ in steps 9, 13, 17, ... The optimality of the strategy was proved for the first 12 steps which is sufficient since thereafter the game is repeating. \\

For $n=7$, with strategy 263265432 the population which were at the beginning in odd numbered boxes are all caught and the repetitive part of the strategy, 6325, produces a relative cat residence probability of $\frac{1}{4}[0, 1, 0, 1, 0, 2, 0]$ in steps 13, 17, 21... But besides its mirror strategy also strategy 225665432$\overline{2563}$ produces the same results: In steps 13, 17, 21, ... the relative cat residence probability is just mirrored: $\frac{1}{4}[0,2, 0, 1, 0, 1, 0]$. The absolute population is reduced by a factor of 16 in every four steps and the average game duration becomes 7373/1792 (4.11440).
\\

For $n=8$ no repetitive strategy pattern was found. The optimal strategy $s_8$ for the first 87 steps is 47527425774224774224774472472552744725527447255274 4725527447255274257752472552742577524~\cite{A3386462}. The resulting average game duration are 4.7496~steps. If any of those digits is changed the game duration will increase.\\

We can observe in Table~\ref{catch_fast} that the game duration is reduced by a factor of 1.210 ($n=5$), 1.274 ($n=6$), 1.342 ($n=7$) and 1.335 ($n=8$) if the faster strategy is applied. An open question is whether this factor will stabilize towards a value around 1.3 or drop maybe even towards 1 when $n$ tends to infinity.

\subsection{Optimal strategy to find cat in a ring}

In this game the two ends of the line of boxes are connected, i.e~if the cat is in box~1 it will move with 50~\% probability to box~2 and with 50~\% probability to box~n. There is no strategy that guarantees to find the cat in a finite number of steps if $n>2$. Nevertheless we can calculate the strategy that finds the cat as fast as possible. Since the boxes are arranged in a ring and all boxes are initially identical, without loss of generality box~1 is opened at the first step.
Table~\ref{catch_fast_in_ring} summarizes the results.

\newpage
\begin{table}[h!]
\begin{center}
\begin{tabular}{|c|c|c|c|}
\hline
 & & \multicolumn{2}{c|}{\textbf{Average duration}}  \\
\textbf{$n$} & \textbf{Strategy}  & \textbf{Game 2.3 (ring)} & \textbf{Game 2.2 (line)} \\
\hline
2 & 11 & 3/2 & 3/2\\
3 & $\overline{1}$ & 7/3 & 5/3\\
4 & $\overline{1}$ & 3.5 & 39/16 (2.4375) \\
5 & $\overline{13524}$ &41/11 (3.$\overline{72}$) & 44/15 (2.9$\overline{3}$) \\
6 & $\overline{14414114}$ & 608/141 (4.31206) & 34165/9984 (3.4220) \\ 
7 & $\overline{1473625}$ & 219/43 (5.09302) & 7373/1792 (4.1144) \\ 
\hline
\end{tabular}
\end{center}
\vspace{-4mm}
\caption{Optimal strategy (or one of the possible optimal strategies) to catch as fast as possible the cat performing a random walk in a ring of $n$ boxes, the average game duration of this strategy and the average game duration of the strategy where box 1 and box n are not connected.}
\label{catch_fast_in_ring}
\end{table}

The case $n=2$ is trivial. In case $n=3$ the player can choose in the first step any box, but in the next steps he must continue opening the same box again and again to achieve the minimum search duration of 7/3.\\

For $n=4$, the player can choose any box in the first step. If the number of the chosen box is odd, then in the following steps they should always stick to odd box numbers.  If the number of the chosen box is even, then in the following steps they should always stick to even box numbers. The game situation is repeating every two steps during which the residence probabilities in each box decrease by a factor of 2. The resulting average search duration is 3.5~steps.\\

For $n=5$, the player can choose box~3 or box~4 in step~2 and then in the following steps they should continue moving two boxes further to the left or to the right. Once the direction is chosen, it must not be changed any more. Hence, the two optimal strategies when starting with box~1 are $\overline{13524}$ and $\overline{14253}$.
The relative residence distribution approaches a steady state with a cycle length of~5. Table~\ref{tab:asymn5_ring} illustrates the cyclic behaviour. The decrease factor of the absolute residence probability between 5 steps approaches 0.18906 which is a root of the polynomial $x^4+42x^3+20x^2-1$. The vector $\Vec{k}$ in Table~\ref{tab:asymn5_ring} is
$$ \Vec{k} = [0.30233, 0.19767, 0.16067, 0.22419, 0.11515].$$

The resulting average search duration is 3.$\overline{72}$~steps.\\

\begin{table}[h!]
\begin{center}
\begin{tabular}{|c|c|c|c|c|c|}
\hline
\textbf{Step $t$} & \textbf{Box 1} & \textbf{Box 2} &\textbf{Box 3} &\textbf{Box 4} &\textbf{Box 5} \\
\hline
$5i+1$ & $k_1$ &$k_2$ & $k_3$ &$k_4$ & $k_5$ \\
$5i+2$ & $k_4$ &$k_5$ & $k_1$ &$k_2$ & $k_3$ \\
$5i+3$ & $k_2$ &$k_3$ & $k_4$ &$k_5$ & $k_1$ \\
$5i+4$ & $k_5$ &$k_1$ & $k_2$ &$k_3$ & $k_4$ \\
$5i$   & $k_3$ &$k_4$ & $k_5$ &$k_1$ & $k_2$ \\
\hline
\end{tabular}
\end{center}
\vspace{-4mm}
\caption{Asymptotic relative residence probability for strategy $\overline{13524}$. Since $k_1$ is the largest number the player should always open the box with the contents~$k_1$.}
\label{tab:asymn5_ring}
\end{table}

For $n=6$, the player must choose box~4 in step~2, but can choose between box~2, box~4 and box~6 in step 3 which all lead to a probability distribution of $\Vec\ell = \frac{1}{48}[\ell_1, 3, \ell_2, 2, \ell_3, 3]$ where one $\ell_i$ is 8 and the other two are 4. Therefore, we get three optimal strategies: $\overline{14252514}$, $\overline{14414114}$ and $\overline{14636314}$. The probability distribution $\Vec\ell$ repeats with a cycle of length 8 and the decrease factor between the 8~steps is $\frac{21}{256}$. The resulting average search duration is 608/141 (4.31206)~steps.\\

For $n=7$, the optimal strategy is similar to the strategy for $n=5$: the player can choose box~4 or box~5 in step~2 and then in the following steps they should continue moving three boxes further to the left or to the right. Once the direction is chosen, it must not be changed any more. Hence, the two optimal strategies when starting with box~1 are $\overline{1473625}$ and $\overline{1526374}$.
The relative residence distribution approaches the steady state [0.21848, 0.12613, 0.11673, 0.15905, 0.07469, 0.18245, 0.12247] which is cyclic in a similar way as the vector~$\Vec{k}$ in Table~\ref{tab:asymn5_ring} but with a cycle length of 7. In each step, the box with the largest residence probability must be opened. The resulting average search duration is 219/43 (5.09302)~steps.\\

\section{Walking grid with exits on left and right side}

If the cat is in the left or right box it will leave the grid with a probability of 50~\%. This is called an escape. If there is no player opening any boxes and catching the cat, the average time it takes the cat starting in box~$i$ to leave the boxes is~\cite{wiki} 
\begin{equation}
    E_{exit}(i,n) = i \cdot (n - i + 1)
\end{equation}

Averaging over all boxes~$i$ yields
\begin{equation}
    E_{exit}(n) = \frac{(n+1)\cdot(n+2)}{6}
\end{equation}

\subsection{Recurring cat residence distribution}

Let $ p_k (i)$ be the probability that the cat is in box $i$ after $k$~steps and $p_k = \sum_{i=1}^n p_k(i)$ the probability that the cat is still in one of the boxes after step~$k$.
If there is no player trying to catch the cat and the initial distribution $p_0$ is not uniform, but follows a sine distribution:
$$ p_0 (i) = A \sin \frac{i \pi}{n+1}$$
with $A=\tan \frac{\pi}{2n+2}$, then the distribution is recurring, only multiplied by a factor of $\cos\frac{\pi}{n+1}$ in each step, because
$$\frac{\sin (x) + \sin(x+h)}{2} = \sin(x + \frac{h}{2})\cdot \cos \frac{h}{2} $$
with $\frac{h}{2} = \frac{\pi}{n+1}$.\\

In this case the probability $p_k$ is 
$$p_k = \cos^k \frac{\pi}{n+1}.$$

The expected number of steps $E_{sin}$ it takes the cat to escape is
\begin{equation}
  E_{sin}(n) = \sum_{i=1}^\infty i \cdot (p_{i-1}-p_i) = \frac{1}{1-\cos\frac{\pi}{n+1}}.  
\end{equation}

As an example we take $n=11$. The cat escapes on average after $E_{sin}(11)=29.35$ steps, but after already 20 steps the probability to be out of the boxes is $\cos^{20} \frac{\pi}{12} = 50$ \%.

\subsection{Opening the boxes at random}

We calculate the escape rate and average game duration if the player opens the boxes in each step at random. Let $E_{ran}(i,n)$ be the escape rate if we have $n$ boxes and the cat starts in box $i$. Since the cat will move from box $i$ either to box $i-1$ or $i+1$ and $1/n$ cats will be removed by the player, the following equations hold for the inner boxes:
$$E_{ran}(i,n) = \frac{n-1}{n}\frac{E_{ran}(i-1,n) + E_{ran}(i+1,n)}{2}.$$
For the first (left) box, the equation reads:
$$E_{ran}(1,n) = \frac{n-1}{n}\frac{1 + E_{ran}(2,n)}{2}.$$
Exploiting the symmetry of the game, we need to solve a linear system of $\left\lceil \frac{n}{2} \right\rceil$~equations. Since the percentage of the cats being caught equals the duration of the game divided by $n$, the following equation holds
\begin{equation}
  \text{escape rate} + \frac{\text{duration}}{n} = 1    
\end{equation}

Table~\ref{open_randomly} summarizes the results. \\

\begin{table}[h!]
\begin{center}
\begin{tabular}{|c|l|l|}
\hline
\textbf{$n$} & \textbf{Escape rate} & \textbf{Average duration} \\
\hline
2 & 1/3 & 1.33 (4/3) \\
3 & 8/21 & 1.86 (13/7) \\
4 & 12/31 & 2.45 (76/31) \\
5 & 0.382 (124/325) & 3.09 (201/65) \\
6 & 0.372  & 3.77 \\
7 & 0.362  & 4.47  \\
8 & 0.351 & 5.19 \\
9 & 0.341 & 5.93 \\
10 & 0.331 & 6.69 \\
20 & 0.262 &  14.77 \\
50 & 0.179 &  41.05 \\
100 & 0.131 & 86.89 \\
1000 & 0.044 & 956.2 \\
\hline
\end{tabular}
\end{center}
\vspace{-4mm}
\caption{Escape rate of a cat performing a random walk through $n$ boxes and average game duration if the player opens the boxes randomly.}
\label{open_randomly}
\end{table}

When simulating this game it can be noticed that the residence distribution is approaching the $\sin \frac{i \pi}{n+1}$ distribution studied in the previous section. However, in this game the decrease factor in each step is $\frac{n-1}{n}\cos\frac{\pi}{n+1}$ since the player removes $\frac{1}{n}$ of the population in each step. Assuming this decrease factor throughout the entire game gives an average game duration of
\begin{equation}
    E_{approx} = \frac{1}{1-\frac{n-1}{n}\cos\frac{\pi}{n+1}}
    \label{eq:e-approx}
\end{equation}
$$$$
which is overestimating the true results, because at the beginning of the game, the cat is not concentrated in the centre of the grid, but as likely resides also in the outer boxes. For $n=1000$, Eq.~\ref{eq:e-approx} gives an average of 995 steps, where the true value is 956.

\subsection{Strategy to minimize the escape rate}

In order to find the strategy that minimizes the probability of escape, a search algorithm was written, that tries all possible strategies up to a length of $maxdepth$ steps, but prunes the branches where "too many cats" have already escaped. For each strategy there is a mirror strategy that yields exactly the same result and therefore we limit the opening of the first box to the left half including the central box in case of $n$ being odd. In all cases up to $n=7$ one or more candidate strategies were found and an obvious repetition pattern identified.\\

Depending on the strategy, the cat residence distribution function starts oscillating after a while with cycles of 4 ($n=6$), 6 ($n=5$),  or 8 ($n=4$). For $n=7$ it is not an exact repetition but it asymptotically approaches a cycle of length 4. Therefore, it is possible to calculate the exact percentage of cats that can escape. The optimal strategies and the escape rates are:\\

\begin{table}[h!]
\begin{center}
\begin{tabular}{|c|c|c|c|}
\hline
\textbf{$n$} & \textbf{Strategy} & \textbf{Escape rate} & \textbf{Average duration} \\
\hline
3 & 22 & 1/3 & 1.33333 \\
4 & $\overline{14414114}$ & 1105/3968 (0.27848) & 2.83543 \\
5 & 144$\overline{141}$ & 9643/37120 (0.25978) & 2.79168 \\
6 & 15261$\overline{2552}$ & 305/1248 (0.24439) & 3.61619 \\ 
7 & 1661$\overline{2266}$ & 183/784 (0.23342) & 4.54974 \\
8 & 177122477$\overline{s_{8r}}$ & 0.22331 & 5.34692 \\
9 & $s_{9}\overline{s_{9r}}$ & 0.21118 & 6.58330 \\
\hline
\end{tabular}
\end{center}
\vspace{-4mm}
\caption{Optimal strategy (or one of the possible optimal strategies in case $n$ is odd) to minimize the escape rate of a cat performing a random walk through $n$ boxes, the resulting escape rate and the average game duration (in case of 8 and 9 boxes at least the first 80 steps of the strategy were proved, the digits of strategies $s_{8r}$, $s_9$ and $s_{9r}$ are given in the text below).}
\label{min_escape_rate}
\end{table}

To prove that the candidate strategy is the optimum, each decision step of the strategy is modified and it is tested if opening a different box at this step will allow to reduce the escape rate. Since we have a tentative minimum given in Table~\ref{min_escape_rate}, any branch of the search tree can be pruned as soon as the escape rate exceeds the presumed minimum. Therefore, for $n < 7$, the proof of the optimal strategy was relatively quick.\\

For $n=3$, opening twice box 2 yields an escape rate of 1/3, in two of three cases the cat will be caught. But also opening box 1 continuously gives the same escape rate. Or opening box 1 for some time and then opening twice box 2.\\

For $n=4$, strategy $\overline{14414114}$ proves to be the only optimum (except for its mirror strategy $\overline{41141441}$). The cat residence distribution function after step 3 is: $\frac{1}{32}[ 3, 2, 3, 1 ]$ and after step 11 it is $\frac{1}{1024}[ 3, 2, 3, 1 ]$. The function has decreased by a factor of 32 in these 8 steps and since the strategy has also a cycle of length 8, this pattern will repeat for ever. Hence we have only to test the first 11 decision steps. And for all 11 steps it turned out that any modification lead to an increase in the escape rate latest at step 16.\\

For $n=5$, there are 4 equivalent strategies (plus their 4 mirror strategies) for the first 9 steps: 14$s_1s_2$141 where $s_{1,2} \in \{25, 41\}$. The residence probability after step 9 always is: $\Vec k_9 = \frac{9}{2560}[ 1, 1, 2, 1, 1 ]$. The proof of optimality of the strategy for the first 9~steps is performed numerically: all strategies which are different lead to an increase in the escape rate latest at step 14. Strategies like 1445 and 1421 have an escape rate larger than 0.26 already after 4 steps.\\

After step 9 there are the two mirror strategies $141$ and $525$ that both return a symmetric residence probability of the form $\frac{9}{20480}\cdot[ 3, 2, 6, 2, 3 ]$ and applied a second time, of the form $\frac{3}{32}\Vec k_9$. Hence we have a repeating cycle for the residence probability of length 6. Checking all possible different strategies up to step 14 where the repetition cycle starts, revealed that all lead to higher escape rates latest after 19 steps. Therefore, the optimal strategies are: 14$s_1s_2$141$s_3 s_4 s_5$... where $s_i \in \{141, 525\}$ for $i>2$.\\

For $n=6$, strategy 15261$\overline{2552}$ proves to be the only optimum (except for its mirror strategy). The cat residence distribution function after step 9 is: $\frac{1}{768}[ 0, 12,  5, 12,  5,  4 ]$ and after step 13 it is $\frac{1}{4096}[ 0, 12,  5, 12,  5,  4 ]$. The function has decreased by a factor of $\frac{16}{3}$ in these 4 steps and this pattern will repeat for ever. Hence we have only to test the first 13 decision steps. And for all 13 steps it turned out that any modification leads to an increase in the escape rate latest at step 29.\\

For $n=7$, strategies 1661$\overline{2266}$ and 1627$\overline{6226}$ are two equivalent strategies. The difference after step~3, where in the first strategy box~6 and in the second box~2 is opened, are in the boxes with odd numbers. In box 2, 4 and 6 the residence probabilities are $\frac{1}{56}[4, 5, 2]$ for both strategies. In the odd boxes it is $\frac{1}{56}[3, 7, 4, 0]$ for the first strategy and $\frac{1}{56}[0, 4, 7, 3]$ for the second, which is just reversed. Since all parts of the odd boxes move in one step from odd to even boxes and vice versa, the "odd" and "even" population never mix. Hence, the equivalent of the first strategy consists of replacing the opening of box~$i$ by box~$7-i$ if the step is odd and if $i$ is even or, if the step is even and if $i$ is odd.\\  

Since both strategies are equivalent, we confine our analysis to strategy 1661$\overline{2266}$. The opening of the boxes repeat with a cycle of length 4 for which we compute the change of the residence probability with the following transition matrix T: 
\[
T = \frac{1}{16}
\begin{bmatrix}
0 & 0 & 0 & 0 & 0 & 0 & 0 & 8 \\
0 & 0 & 0 & 0 & 0 & 0 & 0 & 0 \\
1 & 0 & 3 & 0 & 2 & 0 & 0 & 0 \\
0 & 4 & 0 & 5 & 0 & 2 & 0 & 1 \\
1 & 0 & 4 & 0 & 3 & 0 & 0 & 2 \\
0 & 1 & 0 & 2 & 0 & 1 & 0 & 4 \\
0 & 0 & 1 & 0 & 1 & 0 & 0 & 10 \\
0 & 0 & 0 & 0 & 0 & 0 & 0 & 16 \\
\end{bmatrix}
\]
where $T_{ij}$ defines which fraction of box $i$ is moved to box $j$ in the 4 steps when boxes 2, 2, 6 and 6 are opened. ``Box 8'' counts the fraction of cats that escape. (Note, since all contributions from box 2 are removed when this box is opened in the first step of the repetition cycle, row 2 is empty and could be removed from the calculations. But since the calculations are performed almost instantaneously, we keep it included for better readability.)\\

The matrix T can be diagonalized: 
$$ T=PDP^{-1}$$
where $D$ is the diagonal matrix with the eigenvalues $[0 \ e_1\ e_2\ e_2\ e_1\ 0\ 0\ 1]$ of the matrix $T$. $e_1 = \frac{3-2\sqrt{2}}{16} \approx 0.010723$ and $e_2 = \frac{3+2\sqrt{2}}{16} \approx 0.364277.$\\

The initial residence distribution $p_1$ after opening the boxes 1, 6, 6, 1 becomes:
$$p_1 = \frac{1}{7\cdot 32}[8, 14, 18, 22, 14, 8, 4, 36]$$

The escape rate can be calculated as the last component of the vector $p_\infty$:
$$p_\infty = \lim_{n \to \infty} T^n p_1 = \lim_{n \to \infty} P D^n P^{-1} p_1 $$
Since all eigenvalues except the last one are smaller than 1, $\lim_{n \to \infty} D^n$ is a matrix with zeros except for $D^\infty_{88}=1$. Hence the escape rate becomes
$$ (P^{-1})_{88} \sum_{i=1}^8 P_{i8} \cdot p_{1i} = \frac{183}{784} \approx 0.23341837.$$

Any modification of the conjectured optimal strategies in the first 50 steps leads to an increase in the escape rate latest 11 steps later. A sketch of a proof that this holds for all steps beyond 50 is provided now.\\

The relative residence probability (relative with respect to the remaining fraction of cats still in the boxes) converges towards this repetitive pattern:

\begin{table}[h!]
\begin{center}
\begin{tabular}{|c|c|c|c|c|c|c|c|}
\hline
\textbf{Step $t$} & \textbf{Box 1} & \textbf{Box 2} &\textbf{Box 3} &\textbf{Box 4} &\textbf{Box 5} &\textbf{Box 6} &\textbf{Box 7} \\
\hline
$4i$ & $k_1$ &$k_2$ & $k_3$ &$k_4$ & $k_5$ &$k_6$ & 0 \\
$4i+1$ & 0 & $\ell_1$ &$\ell_2$ & $\ell_3$ &$\ell_4$ & $\ell_5$ &$\ell_6$ \\
$4i+2$ & 0 &$k_6$ & $k_5$ &$k_4$ & $k_3$ &$k_2$ & $k_1$ \\
$4i+3$ & $\ell_6$ & $\ell_5$ &$\ell_4$ & $\ell_3$ &$\ell_2$ & $\ell_1$ & 0 \\
\hline
\end{tabular}
\end{center}
\vspace{-4mm}
\caption{Asymptotic relative residence probability for strategy 1661$\overline{2266}$.}
\label{tab:asymn7}
\end{table}

The factor by which the absolute residence probability decreases between 4 steps approaches $48 - 32 \sqrt{2} = e_2^{-1} \approx 2.745166$. The vector $\Vec{k}$ in Table~\ref{tab:asymn7} is
{\footnotesize $$ \Vec{k} = \frac{1}{1447} [ 893 -563 \sqrt{2}, -177 +319 \sqrt{2}, 660 -233 \sqrt{2}, 147 +152 \sqrt{2}, -233 +330 \sqrt{2}, 157 -5 \sqrt{2}] $$}
or numerically:
$$ \Vec{k} = [0.06690, 0.18945, 0.22840, 0.25015, 0.16150, 0.10361].$$

The vector $\Vec{\ell}$ in Table~\ref{tab:asymn7} is
{\small $$ \Vec{\ell} = \frac{1}{255} [163 -81 \sqrt{2}, 17 +17 \sqrt{2}, 47 +12 \sqrt{2}, 34 +17 \sqrt{2}, -23 +35 \sqrt{2}, 17] $$}
or numerically:
$$ \Vec{\ell} = [0.18999, 0.16095, 0.25086,  0.22761, 0.10391, 0.06667].$$

We analyse first the case of the conjectured optimal strategy after step $4i+1$. Taking the asymptotic relative residence probability given in Table~\ref{tab:asymn7} and applying the strategy $\overline{2662}$ yields an escape rate of $\frac{1}{1785}(103+116\sqrt{2})\approx 0.14961$. Trying out all possible other strategies gives a minimum escape rate of 0.14994 in the following 12 steps.\\

After step 21 of the strategy 1661$\overline{2266}$ for the initial random distribution, we have the distribution $\hat{\ell_j} = \ell_j + \epsilon_j$, $j \in [1, 5]$ in the boxes 2 to 6, where all $\epsilon_j$ are positive except $\epsilon_1$. In box~7 there is $\ell_6$ of which 50~\% escape, identically to the escape in the asymptotic population. Since $\hat{\ell_1}$ is removed by opening box~2 at step $4i+1$ the remaining population is larger than the asymptotic population and hence the escape rate can only increase. This proves that opening box~2 is optimal at step 21.\\

Table~\ref{tab:41} shows how the population propagates from step~$4i+1$ to step~$4i+5$ when applying the strategy 2662. If we look at the contents after 4 more steps of box~3 for example, we can observe that
$$\frac{3\hat{\ell_2}+2\hat{\ell_4}}{16} = \frac{3(\ell_2+\epsilon_2)+2(\ell_4+\epsilon_4)}{16} = \frac{3\ell_2+2\ell_4}{16} + \frac{3\epsilon_2+2\epsilon_4}{16} = \ell_2 \cdot e_2 + \epsilon_{22}.$$
$\epsilon_{22}$ is positive because $\epsilon_2$ and $\epsilon_4$ are positive. The same can be observed in boxes 4, 5 and 6. The contents of box 7 after these 4 steps is $\ell_6 \cdot e_2.$  Hence we have at step~$4i+5$ the same situation as at step~$4i+1$, which means opening box 2 is always optimal for all steps~$4i+i$.\\ 

{\renewcommand{\arraystretch}{1.5}
\begin{table}[h!]
\begin{center}
\begin{tabular}{|c|c|c|c|c|c|c|c|}
\hline
\textbf{Step $t$} & \textbf{Box 1} & \textbf{Box 2} &\textbf{Box 3} &\textbf{Box 4} &\textbf{Box 5} &\textbf{Box 6} &\textbf{Box 7} \\
\hline
$4i+1$ & 0 & $\hat{\ell_1}$ &$\hat{\ell_2}$ & $\hat{\ell_3}$ &$\hat{\ell_4}$ & $\hat{\ell_5}$ &$\ell_6$ \\ \hline
$4i+5$ & 0 & $\frac{4\hat{\ell_3}+\hat{\ell_5}}{16}$ & $\frac{3\hat{\ell_2}+2\hat{\ell_4}}{16}$ & $\frac{5\hat{\ell_3}+2\hat{\ell_5}}{16}$ & $\frac{4\hat{\ell_2}+3\hat{\ell_4}}{16}$ & $\frac{2\hat{\ell_3}+\hat{\ell_5}}{16}$ & $\frac{\hat{\ell_2}+\hat{\ell_4}}{16}$  \\
\hline
\end{tabular}
\end{center}
\vspace{-4mm}
\caption{Propagation for 4 steps of a population that is reached after $4i+1$~steps when applying the strategy 1661$\overline{2266}$ for the initial random distribution.}
\label{tab:41}
\end{table} }

The proof of optimality for the steps~$4i+3$ is identical. Remain the steps $4i$ and $4i+2$ which are symmetric. Taking the asymptotic relative residence probability $\Vec{k}$ and applying the strategy $\overline{2266}$ gives an escape rate of 0.14971. Trying out all possible other strategies gives a minimum escape rate of 0.15014 in the first 12 steps. Applying the strategy 1661$\overline{2266}$ for the initial random distribution, the residence probabilities in box 1, 3 and 5 approach $k_1$, $k_3$ and $k_5$ from the lower side, whereas the residence probabilities in box 2, 4 and 6 approach $k_2$, $k_4$ and $k_6$ from the upper side. This does not allow the same concept of proof as given above for step~$4i+1$. This part remains to be proved.\\

For $n=8$ no short repetitive strategy pattern was found. The optimal strategy for the first 87 steps is 177122477$\overline{s_{8r}}$ with $s_{8r}$ = 2347187237762236818761 (22~digits)~\cite{A3386463}. The escape rate for this strategy is 0.22331 and the average game duration are 5.34692~steps.\\

\vspace{-1mm}
If the strategy $s_{r} = \overline{72347}$ is applied from step~9 onwards, the relative residence probability distribution starts repeating from step 18 onwards with a cycle length of 10. The escape rate for this strategy is only marginally worse: 0.22342. The strategy $s_{r}$ transforms any arbitrary cat residence distribution into a distribution of the form $[5q_1, 0, 15q_1, q_2, 20q_1, q_2, 14q_1, 0]$ after 10 steps. This distribution repeats after another 10 steps, only reduced by a factor of $\frac{1024}{75}$, hence a cycle length of 10 is established.\\

\vspace{-1mm}
For $n=9$ we have a similar set-up as in the case $n=7$, where the "odd" and "even" populations do not mix. The player must open box~1 first (or for the mirror strategy box~9) and in step~2 box~8. In step~3 they have the choice between box~2 and 8. In the first case, the optimal strategy is $s_{9}\overline{s_{9r}}$ with $s_{9}$ = 1829825881238258298723428763 and $\overline{s_{9r}}$ = 9298723458817181238765281318123 876522939298723458297. $s_{9}$ consists of~28 digits and $s_{9r}$ of 52~digits. Like in the case $n=7$ the equivalent strategy is constructed by replacing the opening of box~$i$ by box~$9-i$ if the step is odd and if $i$ is even or, if the step is even and if $i$ is odd. The resulting escape rate is 0.21118 and the average game duration are 6.58330~steps.\\

\vspace{-2mm}
\subsection{Playing strategy "twice left - twice right"}

The strategy $1(n-1)(n-1)1\overline{22(n-1)(n-1)}$ is optimal for $n=7$, but not for $n=8$ and 9 and probably also not for $n>9$. Nevertheless, they are good strategies to contain the escape rate of the cat and they are easy to calculate. Table~\ref{good_strat} gives the escape rates for these strategies, which are an upper limit of the escape rate of the optimal strategy. Obviously, the escape rates are smaller than in cases where the player opens the boxes at random (see Table~\ref{open_randomly}) and therefore it seems plausible that the escape rate also tends towards zero if $n$ tends towards infinity. But this remains to be proved.

\begin{table}[h!]
\begin{center}
\begin{tabular}{|c|l|c|}
\hline
\textbf{$n$} & \textbf{Escape rate} & \textbf{Average duration} \\
\hline
8 & 0.22362 & 5.938 \\
9 & 0.21569 (11/51) &  7.638\\
10 & 0.20849 & 9.644  \\
11 & 0.20274 (133/656) & 11.977 \\
12 & 0.19811 & 14.639 \\
13 & 0.19415 (1999/10296) & 17.629 \\
14 & 0.19073 & 20.947 \\
15 & 0.18778 (10771/57360) & 24.595 \\
16 & 0.18520 & 28.574 \\ 
17 & 0.18292 & 32.884 \\ 
18 & 0.18089 & 37.525 \\ 
19 & 0.17908 & 42.498 \\ 
20 & 0.17745 & 47.803 \\ 
\hline
\end{tabular}
\end{center}
\vspace{-4mm}
\caption{Escape rate and average game duration if a cat is moving randomly between $n$ boxes and the player applies the good but not optimal strategy $1(n-1)(n-1)1\overline{22(n-1)(n-1)}$.}
\label{good_strat}
\end{table}

\section{2d walking grid with no exits}

The analysis will be confined to a 2 x $m$ grid with $m>1$ and to strategies that aim at finding the cat as fast as possible. The boxes in the upper row are numbered from 1 to $m$ and the boxes in the lower row from $m+1$ to $2m$. If the cat has three adjacent boxes, it will choose each box with a probability of 1/3. If it is in one of the four corners, it will move with 50~\% probability to one of the two adjacent boxes.\\

\vspace{-1mm}
In this game, no strategy exists that will guarantee to end the game in a finite number of steps. Every strategy now has 3 mirror strategies due to the vertical and horizontal symmetry. Table~\ref{2d_fast} presents the solutions for $m \le 4$.\\

\begin{table}[h!]
\begin{center}
\begin{tabular}{|c|c|l|}
\hline
\textbf{$m$} & \textbf{Strategy} & \textbf{Average duration} \\
\hline
2 & $\overline{1}$ & 3.5 \\
3 & $\overline{255}$ &  4.11524 \\
4 & 1728$\overline{2772}$ & 5.86092   \\
\hline
\end{tabular}
\end{center}
\vspace{-4mm}
\caption{Optimal strategies and average game duration to catch a cat walking randomly in a 2 x $m$ grid.}
\label{2d_fast}
\end{table}

In case of $m=2$, the cat is confined to a 2 x 2 grid. At odd steps, the player can choose any box. At the next step either the same box or the diagonally opposite box must be opened to minimize the search duration. All these strategies yield an average game duration of 3.5~steps.\\

For $m=3$, the optimal strategy is $\overline{255}$, which leads to a repeating game situation starting after step 3 with a cycle length of 6. The average game duration is 4.11524~steps.\\

For $m=4$, there are two equivalent optimal strategies. Both 1728$\overline{2772}$ and 1771$\overline{7722}$ yield an average game duration of 5.86~steps. The relative residence probability distribution is approaching a steady state with a cycle length of 4 corresponding to the cycle length of the strategies. The proof of optimality was performed until step 36.

\section{2d walking grid with exits in all directions}

\subsection{No opening of boxes}

At a first step we analyse how long it takes on average until the cat leaves a 2~x~$m$~grid when there is no player that opens any boxes. If the cat starts in box $i$ with $i \le \left\lceil \frac{n}{2} \right\rceil$ the average number of steps before the cat leaves the boxes is
\begin{equation}
   E_{2d}(i,m) = 4\; \frac{2a_m-2a_{m-2i}-a_{m-2}+a_{m-2-2i}}{a_{m+1}}
   \label{eq:emi}
\end{equation}

where $a_m$ are the Fibonacci numbers if $m$ is even and the Lucas numbers~\cite{LucasNumber} if $m$ is odd.  Note, to calculate e.g.~$E_{2d}(4,7)$, $a_{-3}$ is required. The terms with negative indices are $a_{-2} = -1$ and $a_{-1} = 1$ for $m$ even and $a_{-3} = -4$ and $a_{-2} = 3$ for $m$ odd.\\

The proof of Eq.~\ref{eq:emi} is straightforward. If the cat is residing in box $i \in [2, \left\lceil \frac{n}{2} \right\rceil -1 ] $, it will move in the next step with the same probability to box $i-1$, to box $i+1$, to the box above or below, which has the same live expectancy as box $i$, and to the outside, which ends the game.\footnote{If the cat resides in box 1 or in a box in the very middle, the proof can be performed in a similar way.} Therefore this equation holds:
\begin{equation}
    E_{2d}(i,m) = \frac{E_{2d}(i-1,m)+E_{2d}(i,m)+E_{2d}(i+1,m)}{4}+1
\end{equation}

we rearrange this equation as 

$$3 E_{2d}(i,m) - E_{2d}(i-1,m)-E_{2d}(i+1,m) = 4$$

Multiplying both sides with $a_{m+1}$ and inserting the right side of Eq.~\ref{eq:emi} yields
\begin{align*}
&24(a_m - a_{m - 2i}) - 12(a_{m - 2} - a_{m - 2 - 2i}) - 8(a_m - a_{m + 2 - 2i}) \\
&\quad + 4(a_{m - 2} - a_{m - 2i}) - 8(a_m - a_{m - 2 - 2i}) + 4(a_{m - 2} - a_{m - 4 - 2i}) = 4a_{m+1}
\end{align*}

Adding the terms gives
\begin{align*}
&8a_m - 4a_{m - 2} - 28a_{m - 2i} +20 a_{m - 2 - 2i} + 8a_{m + 2 - 2i} - 4a_{m - 4 - 2i} = 4a_{m+1}
\end{align*}
Now we exploit that $a_m = 3a_{m - 2}-a_{m - 4}$ and $a_{m+1} = 2a_m-a_{m - 2}$ which leads to 
\begin{align*}
&4a_{m+1} = 4a_{m+1}
\end{align*}
and hence the proof is done.\\

The average over all starting positions gives the average number of steps it takes the cat to exit the grid:

\begin{equation}
    E_{2d}(m) = \frac{4}{a_{m+1}}\left (\frac{2(m-1)}{m}a_m - a_{m-2}\right )
\end{equation}

Table~\ref{2d_dwell_time} presents the average game duration for $m \le 12$. As $m$ tends to infinity, the game duration tends to 4. This is obvious since for an endless grid, the probability that the cat escapes is 1/4 in each step and hence the expected number of steps to escape is 4.

\begin{table}[h!]
\begin{center}
\begin{tabular}{|c|c|c|}
\hline
\textbf{$m$} & \multicolumn{2}{c|}{\textbf{Average duration}} \\
\hline
2 & 2 & 2.000 \\
3 & 52/21  & 2.476 \\
4 & 14/5   & 2.800 \\
5 & 136/45 & 3.022  \\
6 & 124/39 & 3.179   \\
7 & 1084/329 & 3.295   \\
8 & 115/34 & 3.382   \\
9 & 3820/1107 & 3.451  \\
10 & 312/89  & 3.506   \\
11 & 6288/1771 & 3.551   \\
12 &  836/233 & 3.588   \\
\hline
\end{tabular}
\end{center}
\vspace{-4mm}
\caption{Average number of steps a cat takes before escaping when walking randomly in a 2 x $m$ grid.}
\label{2d_dwell_time}
\end{table}

\subsection{Strategy to minimize the escape rate}

For $m=2$ the optimal strategy is the same as when there are no exits: At odd steps, the player can choose any box. At the next step either the same box or the diagonally opposite box must be opened to minimize the escape rate. All these strategies yield an average escape rate of 4/7 and an average game duration of 11/7~steps.\\

It is interesting to note that for any strategy the average game duration is always 1 larger than the escape rate. Let $p_i$ be the probability that the cat was caught at step~$i$. Then the probability that the cat escapes is $\frac{1-p_i}{2}$, the same as the probability that the game does not end with this step because all boxes have exactly a 50~\% chance of escape. Hence, the overall escape rate is 
\begin{equation}
ER = \frac{1-p_1}{2} + \frac{1-p_1}{2}\frac{1-p_2}{2} + \frac{1-p_1}{2}\frac{1-p_2}{2}\frac{1-p_3}{2} ...
\label{eq:escn2}    
\end{equation}

Let  $q_i = p_i +\frac{1-p_i}{2} = \frac{1+p_i}{2}$ be the probability that the game ends during step~$i$, then the average game duration is
\begin{equation}
E_{2d-min} = 1 \times q_1 +  2 \times (1-q_1)q_2 + 3 \times (1-q_1)(1-q_2)q_3 + ...
\label{eq:escn2_av}    
\end{equation}

This sum can be regrouped as
\begin{alignat*}{2}
E_{2d-min} = q_1 + & (1-q_1)q_2 + & (1-q_1)(1-q_2)q_3 + ... \\
 & (1-q_1)q_2 + & (1-q_1)(1-q_2)q_3 + ... \\
 &  & (1-q_1)(1-q_2)q_3 + ... \\
 & & ...
\end{alignat*}
The sum of the first line is 1, the sum of the second line is $1-q_1 = \frac{1-p_1}{2}$, the sum of the third line is $(1-q_1)(1-q_2) = \frac{1-p_1}{2}\frac{1-p_2}{2}$ and so forth. Therefore the difference between the sum in Eq.~\ref{eq:escn2} and the sum in Eq.~\ref{eq:escn2_av} is exactly 1.\\

For $m=3$, the optimal strategy 1$\overline{5522}$ produces an asymptotic residence distribution of $\frac{1}{26}[4, 4, 4, 3, 8, 3]$ which is repeating at every step $4i+1$. The absolute probabilities decrease by a factor of $\frac{256}{7}$ during the 4 steps of the repetition cycle. The escape rate is 0.61797 and the average game duration are 1.86194~steps.\\

For $m=4$, the optimal strategy is to open first a box in a corner, and without loss of generality, we assume box~1 is opened. In step~2 box 7 must be opened and in step~3 the player can choose between box~2 and 7. Both strategies 1728$\overline{2277}$ and 1771$\overline{7227}$ produce the same result: The escape rate is 0.66191 and the average game duration are 2.17652~steps. The results were numerically confirmed until step~45. Table~\ref{min_escape_rate_2d} presents the solutions for $m \le 4$.\\

\begin{table}[h!]
\begin{center}
\begin{tabular}{|c|c|c|c|}
\hline
\textbf{$m$} & \textbf{Strategy} & \textbf{Escape rate} & \textbf{Average duration} \\
\hline
2 & $\overline{1}$ & 4/7 (0.57143) & 11/7 (1.57143)\\
3 & 1$\overline{5522}$ & 0.61797 & 1.86194 \\
4 & 1728$\overline{2277}$ & 0.66191 & 2.17652 \\
\hline
\end{tabular}
\end{center}
\vspace{-4mm}
\caption{Optimal strategies to minimize the escape rate of a cat performing a random walk in a 2 x $m$ grid with exits on all sides, the resulting escape rate and the average game duration.}
\label{min_escape_rate_2d}
\end{table}

Note, the escape rate is increasing with $m$ unlike in the one-dimensional case where the player can sort of control both exits on the left and on the right. 

\section*{Acknowledgements}

The author would like to thank Kester Habermann for performing computer simulations to cross-check the results, Hugo Pförtner for injecting some interesting ideas and mostly Jinyuan Wang for calculating strategies where my Python and Julia programs had not finished after days of branching and pruning the huge serach trees.

\bibliography{references} 

\begin{thebibliography}{1}

\bibitem{Spitzer1964}
Frank Spitzer.
\newblock {\em Principles of Random Walk, 2nd ed.}.
\newblock Springer New York, \url{https://doi.org/10.1007/978-1-4757-4229-9_7}, 1976.

\bibitem{perfect}
Christian Perfect.
\newblock {\em Princess on a Graph}.
\newblock \url{https://www.checkmyworking.com/2011/12/solving-the-princess-on-a-graph-puzzle/}, Retrieved 2025-08-05.

\bibitem{A3386462}
Rüdiger Jehn.
\newblock {\em Sequence A3386462 of the On-Line Encyclopedia of Integer
  Sequences}.
\newblock \url{https://oeis.org/A3386462}, 2025.

\bibitem{wiki}
\newblock \url{https://en.wikipedia.org/wiki/Random_walk}, Retrieved 2025-08-05.

\bibitem{A3386463}
Rüdiger Jehn.
\newblock {\em Sequence A3386463 of the On-Line Encyclopedia of Integer
  Sequences}.
\newblock \url{https://oeis.org/A3386463}, 2025.

\bibitem{LucasNumber}
Eric W Weisstein.
\newblock \url{https://mathworld.wolfram.com/LucasNumber.html}, Retrieved 2025-08-05.

\end{thebibliography}
\bibliographystyle{ieeetr}

\end{document}